\documentclass[11pt]{article}
\usepackage{amsthm}
\usepackage{amssymb}
\usepackage{amsmath}
\usepackage{amsfonts}

\def \remark {\noindent {\bf Remark.}\ \ }
\def \remarks {\noindent {\bf Remarks.}\ \ }

\def \RR {\mathbb R}
\def \EE {\mathbb E}
\def \e {\varepsilon}
\def \eps {\varepsilon}

\def \diam {{\rm diam}}
\def \vol {{\rm Vol}}
\def \Span {{\rm span}}
\def \vrad {{\rm v.rad.}}
\def \Prob {{\rm Prob}}
\def \Med {\text{\sl Med}}

\newtheorem{theorem}{Theorem}[section]
\newtheorem{lemma}[theorem]{Lemma}

\newtheorem{proposition}[theorem]{Proposition}
\newtheorem{corollary}[theorem]{Corollary}

\def \endproof {{\mbox{}\nolinebreak\hfill\rule{2mm}{2mm}\par\medbreak}}

\begin{document}
\title{Small ball probability and Dvoretzky Theorem}
\author{B. Klartag, R. Vershynin}
\date{}
\maketitle

\abstract{Large deviation estimates are by now a standard tool in
the Asymptotic Convex Geometry, contrary to  small deviation
results. In this note we present a novel application of a small
deviations inequality to a problem related to the diameters of
random sections of high dimensional convex bodies. Our results
imply an unexpected distinction between the lower and the upper
inclusions in Dvoretzky Theorem. }

\section{Introduction}

In probability theory, the large deviation theory (or the tail probabilities) and
the small deviation theory (or the small ball probabilities) are in a sense
two complementary directions. The large deviation theory, which is a more
classical direction, seeks to control the probability of deviations of a
random variable $X$ from its mean $M$, i.e. one looks for upper bounds
on $\Prob(|X-M|>t)$. The small deviation theory seeks to control the probability
of $X$ being very small, i.e. it looks for upper bounds on $\Prob(|X|<t)$.
There is a number of excellent texts on large deviations,
see e.g. recent books \cite{DZ} and \cite{dH}.
A recent exposition of the state of the art in the small deviation theory
can be found in \cite{LS}.

A modern powerful approach to large deviations is via the
celebrated concentration of measure phenomenon. One of the early
manifestations of this idea was V. Milman's proof of Dvoretzky
Theorem in the 70s. Recall that Dvoretzky Theorem entails that any
$n$-dimensional convex body has a section of dimension $c \log n$
which is approximately a Euclidean ball. Since Milman's proof, the
concentration of measure philosophy plays a major role in
geometric functional analysis and in many other areas. A recent
book of M. Ledoux \cite{L} gives account on many ramifications of
the method. A standard instance of the concentration of measure
phenomenon is the case of a Lipschitz function on the unit
Euclidean sphere $S^{n-1}$. In view of the geometric applications,
we shall state it for a norm $\|\cdot\|$ on $\RR^n$, or
equivalently for its unit ball, which is a centrally symmetric
convex body $K \subset \RR^n$. Two parameters are responsible for
many geometric properties of the convex body $K$ -- the maximal
and the average value of the norm on the sphere $S^{n-1}$, which
is equipped with the probability rotation invariant measure
$\sigma$:
\begin{equation}                        \label{bM}
b = b(K) = \sup_{x \in S^{n-1}} \|x\|, \ \ \ M = M(K) =
\int_{S^{n-1}} \|x\| \,d\sigma(x).
\end{equation}
The concentration of measure inequality, which appears e.g. in the first pages
of \cite{MS} states that the norm is close to its mean $M$ on most of
the sphere. For any $t > 1$,
\begin{equation}
\sigma \left \{ x \in S^{n-1} :\; \big| \|x\| - M \big| > t M
\right \} < \exp \left(-c t^2 k \right)
 \label{concentration}
\end{equation}
where
$$
k = k(K) = n \left( \frac{M(K)}{b(K)} \right)^2.
$$
Here and thereafter the letters $c,C,c^{\prime},\tilde{c},
c_1,c_2$ etc. denote some positive universal constants, whose
values may be different in various appearances. The symbol
$\asymp$ denotes equivalence of two quantities up to an absolute
constant factor, i.e. $a \asymp b$ if $ca \le b \le Ca$ for some
absolute constants $c, C > 0$.

The concentration of measure inequality can of course be
interpreted as a large deviation inequality for the random
variable $\|x\|$, and the connection to probability theory becomes
even more sound when one recalls an analogous inequality for
Gaussian measures, see \cite{L}. The quantity $k(K)$ plays a
crucial role in high dimensional convex geometry, as it is the
critical dimension in Dvoretzky Theorem.
 We will call this dimension $k(K)$ here the {\em
Dvoretzky dimension}. Milman's proof of Dvoretzky Theorem
\cite{M_dvo} (see also the book \cite[5.8]{MS}) provides accurate
information regarding the dimension of the almost spherical
sections of $K$. Milman's argument shows that if $l < c k(K)$,
then with probability larger than $1 - e^{-c^{\prime} l}$, a
random $l$-dimensional subspace $E \in G_{n,l}$ satisfies
\begin{equation}
{\textstyle \frac{c}{M}} (D^n \cap E)  \  \subset \  K \cap E  \
\subset \ {\textstyle \frac{C}{M}} (D^n \cap E),
\label{inclusion}
\end{equation}
where $D^n$ denotes the unit Euclidean ball in $\RR^n$ and the
randomness is induced by the unique rotation invariant probability
measure on the grassmanian $G_{n,l}$ of $l$-dimensional subspaces
in $\RR^n$.

\medskip The Dvoretzky dimension $k(K)$ was proved in \cite{MS_duke} to be the
exact critical dimension for a random section to satisfy
\eqref{inclusion}, in the following strong sense. If a random
$l$-dimensional subspace $E \in G_{n,l}$ satisfies
(\ref{inclusion}) with probability larger than, say, $1 -
\frac{1}{n}$, then necessarily $l < Ck(K)$. Thus a random section
of dimension $l < c k(K)$ is close to Euclidean with high
probability, and a random section of dimension $l > C k(K)$ is
typically far from Euclidean. These arguments completely clarify
the question of the dimensions in which random sections of a given
convex body are close to Euclidean. Once $b(K)$ and $M(K)$ are
calculated, the behavior of a random section is known. For
instance, Dvoretzky dimension of the cube is $O(\log n)$, while
the cross polytope $K = \{x \in \RR^n :\; \sum |x_i| \leq 1 \}$
has Dvoretzky dimension as large as $O(n)$.

\medskip In this note we investigate Dvoretzky Theorem from a
different direction, which does not involve the standard large
deviations inequality (\ref{concentration}). The second named
author conjectured that a phenomenon similar to the concentration
of measure should also occur for the small ball probability, and
he proved a weaker statement. The conjecture has been recently
proved by R.Latala and K.Oleszkiewitz, using the solution to the
B-conjecture by Cordero, Fradelizi and Maurey \cite{CFM}:

\begin{theorem}[Small ball probability]             \label{lo}
For every $0 < \e < \frac{1}{2}$,
$$
\sigma \left \{ x \in S^{n-1} :\; \| x \| < \e M \right\} < \e^{c
k(K)}
$$
where $c > 0$ is a universal constant.
\end{theorem}

This theorem is related to the small ball probability (as a
direction of the probability theory) in exactly the same way as
the concentration of measure is related to large deviations. Here
we apply Theorem \ref{lo} to study questions arising from
Dvoretzky Theorem. We show that for some purposes, it is possible
to relax the Dvoretzky dimension $k(K)$, replacing it by a
quantity {\em independent of the Lipschitzness} of the norm (which
is quantified by the Lipschitz constant $b(K)$). Precisely, we wish
to replace $k(K)$ by
$$
d(K) = \min \{ -\log \sigma \big \{ x \in S^{n-1} :\; \| x \| \le
{\textstyle \frac{1}{2}} M \big \}, \, n \}.
$$
Selecting $t = \frac{1}{2}$ in the concentration of measure
inequality (\ref{concentration}), we conclude that $d(K)$ must be
at least of the same order as Dvoretzky dimension $k(K)$:
$$
d(K) \ge C k(K).
$$
The small ball Theorem \ref{lo} indeed holds with $d(K)$ (this
is a part of the argument of Latala and Oleszkiewicz, reproduced below).
The resulting inequality can be viewed as Kahane-Khinchine type inequality
for negative exponents:

\begin{proposition}
[Negative moments of a norm]            \label{1-dim Khinchine}
Assume that $0 < l < c d(K)$. Then
$$
cM < \Big( \int_{S^{n-1}} \| x \|^{-l} d \sigma(x)
\Big)^{-\frac{1}{l}} < CM
$$
where $c, C > 0$ are universal constants.
\end{proposition}
For positive exponents, this inequality was proved in \cite{LMS}:
for $0 < k < c k(K)$,
\begin{equation}
c M < \left( \int_{S^{n-1}} \| x \|^k \,d\sigma(x)
\right)^{\frac{1}{k}} < C M.                \label{p_norm intro}
\end{equation}
For negative exponents $-1 < k < 0$, inequality \eqref{p_norm intro}
follows from results of Guedon \cite{G} that generalize
Lovasz-Simonovits inequality \cite{LoSi}.
Proposition \ref{1-dim Khinchine} extends (\ref{p_norm intro}) to
the range $[-c d(K),  c k(K)]$ (which of course includes the range $[-c
k(K), c k(K)]$).

\medskip In Proposition \ref{1-dim Khinchine},
$\|x\|^{-1}$ can be regarded as the radius of the one-dimensional
section of the body $K$. Combining this with the recent inequality
for diameters of sections due to the first named author
\cite{low_M}, we are able to lift the dimension of the section
and thus compute the {\em average diameter of $l$-dimensional
sections} of any symmetric convex body $K$. This formal
``dimension lift'' might be of an independent interest.

\begin{theorem}
[Diameters of random sections] Assume that $0 < l < c d(K)$.
Select a random $l$-dimensional subspace $E \in G_{n,l}$. Then
with probability larger than $1 - e^{-c^{\prime} l}$,
\begin{equation}
{\textstyle K \cap E  \ \subset \ {\textstyle \frac{C}{M}} (D^n
\cap E).} \label{upper_inclusion}
\end{equation}
Furthermore,
$$
\frac{\bar{c}}{M} < \left( \int_{G_{n,l}} \diam(K \cap E)^l \;
d\mu(E) \right)^{\frac{1}{l}} < \frac{\bar{C}}{M}
$$
where $c, c^{\prime}, \bar{c}, C, \bar{C} > 0$ are universal
constants. \label{k-dim Khinchine}
\end{theorem}

The relation between Theorem \ref{k-dim Khinchine} and Dvoretzky
Theorem is clear. We show that for dimensions which may be much
larger than $k(K)$, the upper inclusion in Dvoretzky Theorem
(\ref{inclusion}) holds with high probability. This reveals an
intriguing point in Dvoretzky Theorem. Milman's proof of Dvoretzky
Theorem focuses on the left-most inclusion in (\ref{inclusion}).
Once it is proved that the left-most inclusion in
(\ref{inclusion}) holds with high probability, the right-most
inclusion follows almost automatically in his proof.

\medskip Furthermore, Milman-Schechtman's argument \cite{MS_duke}
implies in fact that the left-most inclusion does not hold (with
large probability) for dimensions larger than the Dvoretzky
dimension. The reason that a random $l$-dimensional section is far
from Euclidean when $l > c k(K)$ is that a typical section does
not contain a sufficiently large Euclidean ball. In comparison to
that, we observe that the upper inclusion in (\ref{inclusion})
 is true in a much wider range
of the dimensions.

\medskip
There are cases, such as the case of the cube, where the Dvoretzky
dimension $k(K)$ is $O(\log n)$, while $d(K)$ is a polynomial in
$n$. Hence, while sections of the cube of dimension $n^c$ are
already contained in the appropriate Euclidean ball (for any fixed
$c < 1$, independent of $n$), only when the dimension is $O(\log n)$ the sections
start to ``fill from inside'', and an isomorphic Euclidean ball is
observed.  The fact that $d(K)$ is typically larger than $k(K)$ is
a bit unexpected. It implies that the correct upper bound for
random sections of a convex body appears sometimes in much larger
dimensions than those in which we have the lower bound.

\medskip In the last decade, diameters of random lower-dimensional
sections of convex bodies attracted a considerable amount of
attention, see in particular \cite{GM1, GM2, GMT}. Theorem
\ref{k-dim Khinchine} is a significant addition to this line of results.
It implies that diameters of random sections are equivalent
for a wide range of dimensions -- starting from dimension one,
when the random diameter simply equals $\frac{1}{M(K)}$, and up to
the critical dimension $d(K)$.

\remark The proof shows that $d(K)$ can be further relaxed in all
our results. For any fixed $u > 1$ it can be replaced by
$$
d_u(K) = \min \{ -\log \sigma \big \{ x \in S^{n-1} :\; \| x \| \le
{\textstyle \frac{1}{u}} M \big \}, \, n \}.
$$

\medskip The rest of the paper is organized as follows.
In Section \ref{s: small ball} we discuss the negative moments of
the norm, proving Proposition \ref{1-dim Khinchine} and Theorem
\ref{lo} by the Latala-Oleszkiewicz's argument. In Section \ref{s:
diameters} we do the lift of dimension and compute the average
diameters of random sections, proving Theorem \ref{k-dim
Khinchine}.

\section{Concentration of measure and the small ball probability}
\label{s: small ball}

We start by proving Proposition \ref{1-dim Khinchine}. It is a
reformulation of the ``small ball probability conjecture'' due to
the second named author. It was recently deduced by R.Latala and
K.Oleszkiewicz from the B-conjecture proved by Cordero, Fradelizi
and Maurey \cite{CFM}. We will reproduce Latala-Oleszkiewicz
argument here. We start with a standard and well-known lemma, on the
close relation between the uniform measure $\sigma$ on the sphere
$S^{n-1}$ and the standard gaussian measure $\gamma$ on $\RR^n$.
For the convenience of the reader, we include its proof.

\begin{lemma}                           \label{transfer}
For every centrally-symmetric convex body $K$,
$$
\frac{1}{2} \sigma(S^{n-1} \cap {\textstyle\frac{1}{2}} K) \le
\gamma(\sqrt{n}K) \le \sigma(S^{n-1} \cap 2K) + e^{-cn}
$$
where $c > 0$ is a universal constant.
\end{lemma}

\proof We will use the two estimates on the Gaussian measure of
the Euclidean ball,
$$
\gamma(2 \sqrt{n} D^n) > \frac{1}{2}, \ \ \
\gamma({\textstyle\frac{1}{2}} \sqrt{n} D^n) < e^{-cn}.
$$
The first estimate is just Chebychev's inequality, and the second
follows from standard large deviation inequalities, e.g. Cramer's
Theorem \cite{V}. Since $K$ is star-shaped,
\begin{align*}
\gamma(\sqrt{n}K)
&\ge \gamma(2 \sqrt{n} D^n \cap \sqrt{n}K) \\
&\ge \gamma(2 \sqrt{n} D^n) \ \sigma_1(2 \sqrt{n} S^{n-1} \cap \sqrt{n}K) \\
\intertext{where $\sigma_1$ denotes the probability rotation
invariant measure on the sphere $2 \sqrt{n} S^{n-1}$} &\ge
\frac{1}{2} \ \sigma(S^{n-1} \cap {\textstyle\frac{1}{2}} K).
\end{align*}
This proves the lower estimate in the Lemma.

For the upper estimate, note that no points of $\sqrt{n}K$ can lie
outside both the ball $\frac{1}{2} \sqrt{n} D^n$ and the positive
cone generated by $\frac{1}{2}S^{n-1} \cap \sqrt{n}K$. Adding the
two measures together, we obtain
\begin{align*}
\gamma(\sqrt{n}K)
&\le \gamma({\textstyle\frac{1}{2}} \sqrt{n} D^n)
  + \sigma_2({\textstyle\frac{1}{2}} \sqrt{n}S^{n-1} \cap \sqrt{n}K) \\
\intertext{where $\sigma_2$ denotes the probability rotation
invariant measure on the sphere $\frac{1}{2} \sqrt{n} S^{n-1}$}
&\le e^{-cn} + \sigma(S^{n-1} \cap 2K).
\end{align*}
This completes the proof.
\endproof

\noindent {\bf Proof of Proposition \ref{1-dim Khinchine}.} As
usual, $K$ will denote the unit ball of the norm $\|\cdot\|$. The
B-conjecture, proved in \cite{CFM}, asserts that the function $t
\mapsto \gamma(e^t K)$ is log-concave. This means that for any $a,
b > 0$ and $0 < \lambda < 1$,
\begin{equation}
\gamma \left( a^{\lambda} b^{1-\lambda} K \right) \geq \gamma
\left( aK \right)^{\lambda} \gamma \left( bK \right)^{1-\lambda}.
\label{concavity}
\end{equation}
Let $\Med = \Med(K)$ be the median of the norm $\|\cdot\|$ on the
unit sphere $S^{n-1}$. By Chebychev's inequality, $\Med \le
2M(K)$. Set $L = \Med \cdot \sqrt{n} K$. According to
Lemma \ref{transfer},
\begin{equation}
\gamma(2L) \ge \frac{1}{2} \sigma(S^{n-1} \cap \Med \cdot K) \ge
\frac{1}{4} \label{gamma_fine1}
\end{equation}
by the definition of the median. On the other hand, again by Lemma
\ref{transfer},
\begin{align}
\gamma({\textstyle\frac{1}{8}} L)
&\le \sigma(S^{n-1} \cap {\textstyle\frac{1}{4}} \Med \cdot K)
   + e^{-cn}  \nonumber \\
&= \sigma(x \in S^{n-1} :\; \|x\| \le {\textstyle\frac{1}{4}} \Med)
   + e^{-cn}                    \label{quarter-median}\\
&\le \sigma(x \in S^{n-1} :\; \|x\| \le {\textstyle\frac{1}{2}}
M(K))
   + e^{-cn}
\le e^{-d(K)} + e^{-c^{\prime}n} < 2 e^{-C d(K)} \nonumber
\end{align}
because $d(K) \leq n$. We may assume that $\eps < e^{-3}$, and
apply (\ref{concavity}) for $a = \eps, b = 2,\lambda =
\frac{3}{\log \frac{1}{\eps}}$. This yields
$$
\gamma( \eps L )^{\frac{3}{\log(1/\e)}}  \gamma(2 L)^{1 -
\frac{3}{\log (1/\e)}} \leq \gamma \left( \eps^{\frac{3}{\log(1/\e)}}
2^{1 - \frac{3}{\log(1/\e)}} L \right)
\leq \gamma({\textstyle\frac{1}{8}} L).
$$
Combining this with (\ref{gamma_fine1})
and (\ref{quarter-median}), we obtain that
$$ \gamma \left( \eps L \right) \leq 8 e^{C^{\prime} d(K) \log
\eps} \leq 8 \eps^{c d(K)} < \left(c^{\prime} \eps \right)^{c
d(K)} $$ and according to Lemma \ref{transfer} we can transfer
this to the spherical measure, obtaining
$$
\sigma(x \in S^{n-1} :\; \|x\| < \e M) < (C\e)^{c d(K)}.
$$
By integration by parts, this yields
$$
\int_{S^{n-1}} \|x/M\|^{-c d(K)/10} d \sigma(x)  \le C,
$$
which implies the left hand side of the inequality
in Proposition \ref{1-dim Khinchine}.
The right hand side follows easily by H\"older's inequality.
\endproof

\medskip By Chebychev's inequality, Proposition \ref{1-dim
Khinchine} yields the desired tail inequality for the small ball
probability:

\begin{corollary}[The small ball probability] \label{cor_lo}
For every $0 < \e < \frac{1}{2}$,
$$
\sigma \left \{ x \in S^{n-1} :\; \| x \| < \e M \right\} < \e^{c
d(K)} < \e^{c^{\prime} k(K)}.
$$
where $c, c^{\prime} > 0$ are universal constants.
\end{corollary}

Theorem \ref{lo} is contained in Corollary \ref{cor_lo}.
Let us give some interpretation of the expression in Proposition
\ref{1-dim Khinchine}. For a subspace $E \subset \RR^n$, let $S(E)
= S^{n-1} \cap E$ and $\sigma_E$ be the unique rotation invariant
probability measure on the sphere $S(E)$. We will use the fact
that $\vol(K) = \vol(D^n) \int_{S^{n-1}} \|x\|^{-n} \,d\sigma(x)$.
The {\em volume radius} of a $k$-dimensional set $T$ is defined as
$$
\vrad(T) = \left( \frac{\vol(T)}{\vol(D^k)} \right)^{\frac{1}{k}}.
$$
Thus
$$
\vrad(K) = \left( \int_{S^{n-1}} \|x\|^{-n} \;d\sigma(x) \right)^{1/n}.
$$
By the rotation invariance of all the measures (as in \cite{low_M}), we
conclude that
\begin{align}
\int_{S^{n-1}} \|x\|^{-k} \,d\sigma(x)
&= \int_{G_{n,k}} \int_{S(E)} \|x\|^{-k} \,d\sigma_E(x) \,d\mu(E) \nonumber\\
&= \int_{G_{n,k}} \vrad(K \cap E)^k \,d\mu(E).    \label{vrad}
\end{align}
Thus Proposition \ref{1-dim Khinchine} asymptotically computes the
average volume radius of random sections. This fits perfectly to
the estimates for diameters of sections in \cite{low_M}, to be
applied next.

\section{Diameters of random sections}
\label{s: diameters}

In this section we prove the main result of the paper, Theorem
\ref{k-dim Khinchine}. We regard $\|x\|^{-1}$ as the radius of the
one-dimensional section spanned by $x$; thus Proposition
\ref{1-dim Khinchine} is an asymptotically sharp bound on the
diameters of random one-dimensional sections. Theorem \ref{k-dim
Khinchine} extends this bound to $k$-simensional sections, for all
$k$ up to the critical dimension $d(K)$. We start with a ``lift of
dimension'', which is a variation of the ``low $M$ estimate'',
Proposition 3.9 of \cite{low_M} (the case $\lambda = \frac{1}{2}$
there). The difference to that Proposition, is that here we
estimate the $L_k$ norm, rather than just the tail probability as
in \cite{low_M}.

\begin{proposition}
[Dimension lift for diameters of sections]      \label{dimension
lift} Let $1 \le k_0 < n$. Then for any integer $k < k_0/4$,
$$
\left( \int_{G_{n,k}} \diam(K \cap E)^k \; d\mu(E)
\right)^{\frac{1}{k}} \le C M(K) \left( \int_{S^{n-1}}
\|x\|^{-k_0} \; d\sigma(x) \right)^{\frac{2}{k_0}}.
$$
\end{proposition}

\remarks

{\bf 1.}
 Theorem \ref{k-dim Khinchine} follows immediately from
Proposition \ref{1-dim Khinchine} and Proposition \ref{dimension
lift}. Indeed, the right hand side in Proposition \ref{dimension
lift} is bounded by
$$ CM(K) \left(\frac{C}{M(K)} \right)^2 \le \frac{C}{M(K)}. $$

{\bf 2.} Note the special normalization in Proposition
\ref{dimension lift} (compare with Proposition 3.9 in
\cite{low_M}). Actually, as follows from the proof, for any
$\lambda > 0$, the right hand side in Proposition \ref{dimension
lift} may be replaced with
$$ C(\lambda) M(K)^{\lambda} \left( \int_{S^{n-1}}
\|x\|^{-k_0} \; d\sigma(x) \right)^{\frac{1 + \lambda}{k_0}}. $$
in the price of increasing $k_0$ (just replace in the proof
Cauchy-Schwartz with the appropriate H\"older inequality).
However, we cannot have $\lambda = 0$, as is demonstrated by the
example of $K = \RR^{n-1}$. The average diameter of
one-dimensional sections of $K$ is zero, while the diameter of any
section of dimension larger than one is infinity. Thus there can
not be any formal dimension lift unless one has an extra factor
which must be infinity for ``flat'' bodies. Such a factor is $M(K)
= \int_{S^{n-1}} \|x\| \,d\sigma(x)$.

\medskip In order to prove Theorem \ref{1-dim Khinchine} we
need a standard lemma that claims that the average norm $M$ is
stable under the operation of passing to a random subspace. We are
unaware of a reference for the exact statement we need (a similar
result appears e.g. in Lemma 6.6 of \cite{M_iso}), so a proof is
provided. The average norm on a subspace $E \in G_{n,k}$ is denoted
by $M_E = \int_{S(E)} \|x\| \,d\sigma_E(x)$.

\begin{lemma}
For every norm $\| \cdot \|$ on $\RR^n$ and every integer $0 < k < n$,
\begin{equation}            \label{raz_lemma eq}
c M <
\left( \int_{G_{n,k}} (M_E)^{2k} d\mu(E) \right)^{\frac{1}{2k}}
< C M
\end{equation}
where $c, C > 0$ are universal constants.    \label{raz_lemma}
\end{lemma}

\proof The first inequality in \eqref{raz_lemma eq} follows easily
from H\"older's inequality. In the proof of the second inequality,
we will use a variant of Raz's argument (see \cite{R, MW}). We
normalize so that $M = 1$. Let $X_1,..,X_k$ be $k$ independent
random vectors, distributed uniformly on $S^{n-1}$. It is well
known that a norm of a random vector on the sphere has a
subgaussian tail (e.g. \cite[3.1]{LMS}):
$$
\EE \exp \left (s \|X_i\| \right) < \exp \left(c s^2 \right)
\ \ \ \text{for all $i$ and all $s > 1$}
$$
which by independence implies
$$
\EE \exp \left(s \cdot \frac{1}{k} \sum_{i=1}^k \| X_i \| \right)
< \exp \left( \frac{C s^2}{k} \right)
\ \ \ \text{for $s > 1$.}
$$
Using Chebychev's inequality and optimizing over $s$ (e.g.
\cite[7.4]{MS}), we get
\begin{equation}
\Prob \left \{ \frac{1}{k} \sum_{i=1}^k \| X_i \| > C t \right \}
< \exp \left(-t^2 k \right)
\ \ \ \text{for $t > 1$.}           \label{basic_A}
\end{equation}
Let $E$ be the linear span of $X_1,..,X_k$. Then $E$ is distributed
uniformly in $G_{n,k}$ (up to an event of measure zero). Since for
any two events one has $\Prob(A) \leq \frac{\Prob(B)}{\Prob(B|A)}$,
we conclude that
\begin{equation}
\Prob \left \{M_E > 2ct \right \}
\leq \frac{\Prob \left\{ \frac{1}{k} \sum_{i=1}^k \| X_i \| > c t \right \}}
       {\Prob \left\{\frac{1}{k}\sum_{i=1}^k \| X_i \| > c t
         \, \big| \, M_E > 2 ct \right \}}.
\label{condition}
\end{equation}
The enumerator in \eqref{condition} is bounded by \eqref{basic_A}.
To bound below the denominator note that
$\| X_i \| < C \sqrt{k} M_E$ pointwise for all $i$.
This is a consequence of a simple comparison inequality for the
Gaussian analogs of $M$, $M_E$ (see e.g. \cite[5.9]{MS}).
Note that each $X_i$ is distributed uniformly in $S(E)$.
For a fixed $E \in G_{n,k}$, we can estimate
$P_E = \left\{ \frac{1}{k} \sum_{i=1}^k \| X_i \| > \frac{M_E}{2}
\,\big|\,  \Span \{ X_1,..,X_k \} = E \right\}$ via Chebychev's
inequality as
$$
M_E = \EE  \left( \frac{1}{k} \sum_{i=1}^k \| X_i \|
\,\Big|\,  \Span \{ X_1,..,X_k \} = E \right)
\le C \sqrt{k} M_E P_E + \frac{M_E}{2} (1 - P_E).
$$
Hence $P_E \ge \frac{\bar{c}}{\sqrt{k}}$ for every $E \in G_{n,k}$.
Thus
\begin{align*}
\text{denominator in \eqref{condition}}
&\ge \Prob \left \{ \frac{1}{k} \sum_{i=1}^k \| X_i \|
     > \frac{M_E}{2} \,\Big|\, M_E > 2 ct \right \} \\
&\ge \min_{E \in G_{n,k}} P_E
\ge \frac{\bar{c}}{\sqrt{k}}.
\end{align*}
Combining with (\ref{basic_A}) and (\ref{condition}) we get
$$
\Prob \left \{M_E > 2 ct \right\}
< c \sqrt{k} e^{-t^2k}
< e^{-C t^2 k}
\ \ \ \text{for $t > 1$.}
$$
Integrating by parts we get the desired estimate.
\endproof

\emph{Proof of Proposition \ref{dimension lift}.} By H\"older
inequality, the right hand side increases with $k_0$, hence we may
assume that $k_0 = 4k$. We shall rely on the main result of
\cite{low_M}, which claims that for any centrally-symmetric convex
body $T \subset \RR^n$, and every integers $0 < k \leq l < n$,
\begin{equation}                        \label{boaz}
\vrad(T) > C \left( \int_{G_{n,k}} \vrad(T \cap E)^l \, \diam(T
\cap E)^{n-l}
  \;d\mu(E) \right)^{1/n}.
\end{equation}
We are going to apply \eqref{boaz} to $T = K \cap E$, for
subspaces $E \in G_{n,k_0}$. Denote by $G_{E,k}$ the Grassmanian
of all $k$-dimensional subspaces of $E$, equipped with the unique
probability rotational invariant measure. Then by \eqref{boaz}, \eqref{vrad}
and the rotational invariance of all measures,
\begin{eqnarray*}
\lefteqn{ \int_{E \in G_{n,k}} \vrad(K \cap E)^{2k} \, \diam(K
\cap E)^{2k} \;d\mu(E) } \\ & = &
 \int_{E \in G_{n,4k}} \int_{F \in G_{E,k}} \vrad(K \cap F)^{2k} \, \diam(K \cap F)^{2k}
   \; d\mu_E(F) d\mu(E) \\ & \leq & C^{k_0} \int_{E \in G_{n,k_0}}
  \vrad (K \cap E)^{k_0} d\mu(E) = C^{k_0} \int_{S^{n-1}} \| x \|^{-k_0}
  d\sigma(x).
  \label{diam_vr}
  \end{eqnarray*}
Also, by Cauchy-Schwartz inequality,
\begin{eqnarray*}
\lefteqn{ \left( \int_{E \in G_{n,k}} \diam(K \cap E)^{k}
\;d\mu(E) \right)^{\frac{1}{k}} }
\\ & \leq &
 \left( \int_{E \in G_{n,k}} \vrad(K \cap E)^{2k} \, \diam(K \cap E)^{2k}
  \;d\mu(E) \right)^{\frac{1}{2k}} \left( \int_{E \in G_{n,k}} \frac{1}{\vrad(K \cap
  E)^{2k}} \;d\mu(E) \right)^{\frac{1}{2k}}.
  \end{eqnarray*}
We will use the standard inequality $\frac{1}{v.rad.(K \cap E)}
\leq M_E$, which follows directly from H\"older inequality (e.g. \cite{MS}).
Then,
\begin{eqnarray*}
\lefteqn{ \left( \int_{E \in G_{n,k}} \diam(K \cap E)^{k}
\;d\mu(E) \right)^{\frac{1}{k}} } \\ & \leq & C \left(
\int_{S^{n-1}} \| x \|^{-k_0} d\sigma(x) \right)^{\frac{2}{k_0}}
\left( \int_{E \in G_{n,k}} ( M_E )^{2k} \;d\mu(E)
\right)^{\frac{1}{2k}}
  \end{eqnarray*}
and the proposition follows by Lemma \ref{raz_lemma}.
  \endproof

\end{document}